\documentclass[11pt]{article}
\usepackage{amsthm}
\usepackage{amsmath,amssymb}
\usepackage{colordvi}
\usepackage[all]{xy}

\newtheorem{thm}{Theorem}
\newtheorem{prop}[thm]{Proposition}
\newtheorem{cor}[thm]{Corollary}

\newcommand{\C}{{\mathbb C}} 
 
\newcommand{\Z}{{\mathbb Z}} 
\newcommand{\N}{{\mathbb N}} 
\newcommand{\Q}{{\mathbb Q}} 
\newcommand{\PP}{{\mathbb P}}

\newcommand{\Ae}{{\mathcal A}}
\newcommand{\Ce}{{\mathcal C}}
\newcommand{\De}{{\mathcal D}}
\newcommand{\Ge}{{\mathcal G}}
\newcommand{\Je}{{\mathcal J}}
\newcommand{\Le}{{\mathcal L}}
\newcommand{\Me}{{\mathcal M}}

\newcommand{\We}{{\mathcal W}}

\newcommand{\Pg}{{\mathfrak P}}

\newcommand{\Aut}{{\rm Aut}}

\newcommand{\Gal}{{\rm Gal}}

\renewcommand{\div}{{\rm div}}
\newcommand{\et}{^{\ast}}

\newcommand{\Jac}{\Je\!ac}

\newcommand{\bibname}{}

\title{The Weierstrass subgroup of a curve\\has maximal rank.}
\author{Martine Girard, David R. Kohel, and Christophe Ritzenthaler}

\date{}
\begin{document}

\maketitle

\begin{abstract}
We show that the Weierstrass points of the generic curve of genus~$g$ 
over an algebraically closed field of characteristic~$0$ generate a 
group of maximal rank in the Jacobian.
\end{abstract}

The Weierstrass points are a set of distinguished points on curves,
which are geometrically intrinsic.
In particular, the group these points generate in the Jacobian is a
geometric invariant of the curve. A natural question is to determine
the structure of this group. For some particular curves with large
automorphisms groups (for instance, Fermat curves~\cite{Rohrlich77}),
these groups have been found to be torsion. 
The first author provided the first examples where this group has
positive rank (\cite{G:JNT}, \cite{G:Acta}) and obtained a lower
bound of 11 on the rank of the generic genus 3 curve. 
The motivation of this paper was to bridge the gap between this bound and 
the expected bound of 23 -- meaning that there are no relations
between the Weierstrass points on the generic genus 3 curve. 
The result we obtain is valid for generic curves of any genus. More
precisely, let the {\it Weierstrass subgroup} of a curve $\Ce$ be the
group generated by the Weierstrass points in the Jacobian of the curve
$\Ce$. We show that  

\begin{thm}\label{Main}
  The Weierstrass subgroup of the generic curve of genus $g \geq 3$ is 
  isomorphic to $\Z^{g(g^2-1)-1}$.
\end{thm}

\noindent
As a consequence of this theorem, we deduce the following corollaries.

\begin{cor}
  For any number field $K$, the group generated by the Weierstrass 
  points of a curve over $K$ in its Jacobian is isomorphic to 
  $\Z^{g(g^2-1)-1}$, outside of a set of curves whose moduli lie in 
  a thin set in $\Me_g(K)$.
\end{cor}

\begin{cor}
  For each $g \le 13$, the curves of genus $g$ over $\Q$ for which 
  the group generated by the Weierstrass points in its Jacobian is 
  isomorphic to $\Z^{g(g^2-1)-1}$, determine a Zariski dense set of 
  moduli in $\Me_g$.
\end{cor}

We start by recalling some basic definitions and properties of
Weierstrass points, then some results concerning the behaviour of
Weierstrass points under specialisation. We then describe the 
fundamental tools in our study, which are the natural Galois module
structure of the subgroup of divisors with support on the 
Weierstrass points and the geometric characterisation of the Galois 
group.  Using the specialisation of this Galois module in families,
we obtain the main result.

\section{The Weierstrass subgroup of a curve.}

We recall in this section the definition and some properties of
Weierstrass points (see~\cite{HS00} exercise A.4.14).

Let $\Ce/K$ be a smooth projective curve of genus $g \geq 2$ over a field 
$K$ of characteristic $0$, and let $P$ be any point on $\Ce$.  We say that 
$P$ is a {\em Weierstrass point} if and only if there exists a differential 
form $ \omega \in H^{\circ}(C,\Omega_C)$, such that ${\rm ord}_P(\omega) 
\geq g$. Let $\We$ be the set of Weierstrass points on $\Ce$.

Alternatively we may characterise the Weierstrass points in terms of 
the dimensions of Riemann-Roch spaces.  For any divisor $D$ on $\Ce$, 
let $\Le(D)$ be the Riemann-Roch space 
$$
\{f \in K(\Ce)\et \,:\, \div(f)+D \geq 0\} \cup \{0\}
$$ 
and let $\ell(D)$ be its dimension.  

\begin{prop}
  A point $P$ is Weierstrass if and only if $\ell(gP)\geq 2$.
\end{prop}

\noindent
The gap sequence associated to a point $P$ is defined to be the set 
$$
G(P) = \{ n \in \N \,:\, \ell(nP)=\ell((n-1)P) \}.
$$  
We can define the {\em weight} of a point to be 
$$
w(P) = \Big(\sum_{n \in G(P)} n\Big) - g(g+1)/2.
$$
Positive weight provides yet another characterisation of Weierstrass 
points and we have a formula for the number of Weierstrass points,
counted with multiplicities equal to their weights.

\begin{prop}
  A point $P$ is a Weierstrass point if and only if $w(P) \geq 1$, 
  and $\sum w(P)P$ belongs to the complete linear system
  $$
  \left|\frac{g(g+1)}{2}K_C\right|.
  $$
  In particular, the sum of the weights of all Weierstrass points 
  is $g(g^2-1)$.
\end{prop}

We define the {\em Weierstrass subgroup} to be the group generated 
by the differences of Weierstrass points in the Jacobian of the curve
identified with ${\rm Pic}^{\circ}(\Ce)$.

\section{Specialisation of Weierstrass points.}

In this section, we recall two theorems on the behaviour of 
Weierstrass points under specialisation.   
For a family $\Ce \rightarrow X$ of smooth projective curves of 
genus~$g$ over an irreducible base $X$.  We define $W_{\eta}$ and 
$W_s$ to be the group generated by the Weierstrass points in the 
generic fibre $\Ce_\eta$ and a special fibre $\Ce_{s}$, respectively.

\begin{thm}
  \label{thm_inject}
  The group of Weierstrass points form an algebraic family 
  such that $W_{\eta}$ surjects on $W_{s}$ and is injective 
  on torsion.
\end{thm}

\begin{proof}
  For the first part, see Hubbard~\cite{Hubbard76} or 
  Laksov-Thorup~\cite{LaksovThorup94}, the second part is classic 
  (see~\cite{HS00} Theorem C.1.4).   
\end{proof}

We furthermore need the following theorem of N{\'e}ron which provides 
constraints on the set of rational points for which the generic rank 
can decrease under specialisation.

\begin{thm}[N{\'e}ron~\cite{Neron52}, see Serre~\cite{Serre89}, p.152]
  \label{thm_neron}
  Let $\pi: \Ae \rightarrow X$ be a family of abelian varieties 
  over a number field $K$.  Then the group $\Ae(K(X))$ is finitely 
  generated, and the set 
  $$
  \{ P \in X(K) \;|\;  \Ae(K(X)) \rightarrow \Ae_P(K) 
  \mbox{ is not injective} \}
  $$
  is a thin set.
\end{thm}

\noindent
{\bf N.B} This form of the theorem appears in Serre~\cite{Serre89} 
with $X$ any open subvariety of $\PP^n$, but the proof holds more 
generally.  
However, for $X$ of dimension 1 or of general type, the full set 
$X(K)$ of rational points may be thin for any number field $K$.  

\section{The generic Galois group.}

Let $\Me_g$ be the moduli space of curves of genus $g$, let $\Ce_g \to \Me_g$ 
be the universal curve over $\Me_g$ (i.e. the moduli space of pointed curves), 
and let $\pi: \We_g \to \Me_g$ be the restriction to the locus of Weierstrass 
points.  Eisenbud and Harris~\cite{EisenbudHarris87a} studied the geometric 
monodromy group of this cover, which can be identified 
(see Harris~\cite{Harris79}) with the geometric Galois group of $\pi$, i.e. 
the group of automorphisms of the normal closure of $\C(\We_g)/\C(\Me_g)$.  
They proved that the monodromy group is as large as possible. 

\begin{thm}[Eisenbud-Harris~\cite{EisenbudHarris87a}]
  \label{thm_monodromy}
  The monodromy group of $\pi$ equals the full symmetric group $S_{g(g^2-1)}$ 
  acting on the $g(g^2-1)$ generic Weierstrass points.
\end{thm}

Since $\Me_g$ can be defined over $\Q$, and the geometric Galois group 
is maximal, we conclude that the Galois group of the normal closure of 
$\Q(\We_g)/\Q(\Me_g)$ must also be the full symmetric group.  We apply 
this theorem to the Weierstrass subgroup of the Jacobian, as a Galois 
module, in order to prove Theorem~\ref{Main}. 

\section{Galois module structure.}

Let $\Ce/K$ be a curve of genus $g$, and let $\We$ be its set of 
Weierstrass points in $\Ce(\bar{K})$.  Then the absolute Galois 
group $\Ge = \Gal(\bar{K}/K)$ acts on the set $\We$.  Thus the 
Weierstrass divisor group,
$$
V = \bigoplus_{P \in \Pg}\ \Z. P
$$
is equipped with a natural $\Z[\Ge]$-module structure, which acts 
through $\Z[G]$, where $G$ is the image of $\Ge$ in $\Aut(\We) = 
S_{g(g^2-1)}$ acting as permutations of $\We$.  

\begin{thm}
  \label{thm_free_or_tors}
  The Weierstrass subgroup of the Jacobian $\Jac(\Ce_g)$ of the 
  generic curve $\Ce_g$ is either a free group of rank $g(g^2-1)-1$
  or a torsion subgroup.
\end{thm}

\begin{proof}
Both the ``degree zero'' submodule $V^\circ$ of $V$, generated by 
differences of Weierstrass points, and its submodule $P$ of principal 
divisors with support in $\We$ are $\Z[\Ge]$-submodules.
From Theorem~\ref{thm_monodromy}, we know that $\We$ consists of one 
orbit of $\Ge$, which acts through the full symmetric group $S_n$, 
where $n = g(g^2-1) = |\We|$.  Thus $P_\Q = P\otimes_\Z\Q$ is a 
$\Q[S_n]$-submodule of $V^\circ_\Q = V^\circ\otimes_\Z\Q$.  
Since $V^\circ_\Q$ is simple as a $\Q[S_n]$-module, it follows that 
$P_\Q$ is either trivial or equal to $V^\circ_\Q$, and the theorem 
follows.
\end{proof}

\section{Weierstrass subgroups of cyclic covers of $\PP^1$.}

In this section, we find curves of any genus $g$ such that the subgroup 
generated by the difference of two Weierstrass points $P$ and $Q$ 
has odd order.  Comparing this with the Weierstrass subgroup of a 
hyperelliptic curve, we establish that the generic Weierstrass subgroup 
can not be a torsion subgroup.  First, we state the classical result 
for the Weierstrass subgroup of a hyperelliptic curve.

\begin{prop}
  \label{prop_hyperelliptic}
  The Weierstrass points of an hyperelliptic curve of genus $g$ generate
  the $2$-torsion subgroup of its Jacobian.
\end{prop}

In any genus, there exists cyclic trigonal covers of the projective line.
Such covers can be described (see~\cite{Coppens85} and \cite{Kato80}) by
the plane model 
$$
C : y^3 = \prod_{i=1}^s (x-\alpha_i) \prod_{j=1}^t (x-\beta_j)^2
$$
where $\alpha_i$ and $\beta_j$ are distinct complex numbers and 
$s$ and $t$ satisfy $s+2 t \equiv 0 \bmod 3$ and $t<s$. 
The genus of $C$ is then equal to $g=s+t-2$.
 
\begin{prop}
  \label{prop_trigonal}
  There exists a curve of genus $g$ with two Weierstrass points whose
  difference is a point of order $3$ in the Jacobian.
\end{prop}

\begin{proof}
  We take $C$ a trigonal curve as defined above, of genus $g > 2$, 
  with $t$ in $\{0,1,2\}$ such that $t \equiv -g + 1 \pmod{3}$, 
  and with $s = g - t + 2 \ge 2$.  Then there exist two nonsingular 
  points $P_1 = (\alpha_1,0)$ and $P_2 = (\alpha_2,0)$. 
  The functions $f = (x-\alpha_1)/(x-\alpha_2)$ and $1/f$ are 
  respectively in $\Le(3 P_2)$ and in $\Le(3 P_1)$, and thus 
  the points $P_1,P_2$ are Weierstrass points. 
  Moreover, since $\div(f) = 3 (P_1 - P_2)$ it follows that 
  $P_1-P_2$ is a $3$-torsion point in the Weierstrass subgroup of 
  the Jacobian of $C$. 
\end{proof}

\section{Proof of the main theorem.}

We are now in a position to prove:

\setcounter{thm}{0}
\begin{thm}
  The Weierstrass subgroup of the generic curve of genus $g \geq 3$ is 
  isomorphic to $\Z^{g(g^2-1)-1}$.
\end{thm}
\begin{proof}
  By Theorem~\ref{thm_free_or_tors}, the generic Weierstrass subgroup 
  is either a free abelian group or is purely torsion.
  In the latter case, Theorem~\ref{thm_inject} implies that the 
  generic Weierstrass subgroup is isomorphic with the Weierstrass subgroup 
  of every special curve. 
  We first consider the moduli space $T = \Me_g^{(m)}$ with $m$-level
  structure.  For $m \ge 3$, the space $T$ is a fine moduli 
  space, with universal cover $\De_g \rightarrow T$ 
  such that each fibre $C_t$ is a curve of genus $g$ whose isomorphism 
  class determines the moduli point $\pi(t)$ on $\Me_g$~(see 
  e.g.~\cite{GeemOort00}). 
  On the other hand, the finite covers $\De_g \rightarrow \Ce_g$ 
  and $T \rightarrow \Me_g$
  determine a birational morphism of $\De_g$ to the fibre 
  product $\Ce_g \times_{\Me_g} T$, by which we may identify 
  the generic Weierstrass subgroup $W_\eta$ with the generic Weierstrass 
  subgroup over $T$ (since then $\De_g/\Q(T)$ is isomorphic to 
  $\Ce_g/\Q(T)$).
  By specialising to a hyperelliptic curve $C_t$, 
  Proposition~\ref{prop_hyperelliptic} implies that, if torsion, the 
  generic Weierstrass subgroup must equal the $2$-torsion subgroup.  
  This contradicts the result of Proposition~\ref{prop_trigonal} which 
  implies that it must surject on a subgroup of order $3$.  
  We conclude that $W_{\eta}$ contains a point of infinite order, 
  and thus is free of maximal rank.
\end{proof}

In the following corollaries we apply the main theorem to the 
specialisations of the Weierstrass subgroup to a number field $K$.  
Since the generic Weierstrass subgroup is defined only over the 
splitting field of $\Q(\We_g)/\Q(\Me_g)$, we apply Weil restriction 
to the generic Jacobian to find a family of abelian varieties over 
$\Me_g$ to which we can apply N\'eron's specialisation theorem.  
As above, we let $\Jac(\Ce_g)/\Q(\Me_g)$ be the Jacobian of the generic 
fibre of $\Ce_g \rightarrow \Me_g$. Set $n = g(g^2-1)$, let $\Me_g^\circ$ 
be the open subvariety of $\Me_g$ on which $\Ce_g$ has $n$ distinct 
Weierstrass points, and let $\We_{g,n}$ be the moduli space of curves 
with $n$ ordered Weierstrass points, as a finite cover of $\Me_g^\circ$.  
Now let $A/\Q(\Me_g)$ be the Weil restriction of $\Jac(\Ce_g)$ with 
respect to the extension $\Q(\We_{g,n})/\Q(\Me_g)$.  
Then $A(\Q(\Me_g))$ is naturally isomorphic to $\Jac(\Ce_g)(\Q(\We_{g,n}))$, 
which contains the generic Weierstrass subgroup.  
By restricting $\Me_g^\circ$ further to an open subvariety, we obtain 
a family $\Ae \rightarrow \Me_g^\circ$ with generic fibre $A/\Q(\Me_g)$.

\begin{cor}
  For any number field $K$, the group generated by the Weierstrass 
  points of a curve over $K$ in its Jacobian is isomorphic to 
  $\Z^{g(g^2-1)-1}$, outside of a set of curves whose moduli lie in 
  a thin set in $\Me_g(K)$.
\end{cor}

\begin{proof}
We apply N{\'e}ron's specialisation theorem (Theorem~\ref{thm_neron}) 
to the family $\Ae \rightarrow \Me_g^\circ$ above.  For every point $P$ 
in $\Me_g^\circ(K)$ outside of a thin set, the group $\Ae(K(\Me_g)) = 
A(K(\Me_g))$ --- and the Weierstrass subgroup in particular --- 
specialises injectively to $\Ae_P(K)$.  
Since $\Me_g^\circ$ is open in $\Me_g$, the same holds for every $P$ 
outside of a thin set in $\Me_g(K)$.
\end{proof}

For $g \le 6$ the moduli space $\Me_g$ is rational, and the complement 
of a thin set in $\Me_g(\Q)$ provides a dense set in $\Me_g$ consisting  
of moduli of curves whose Weierstrass subgroup has maximal rank.
More generally, this latter property holds in $\Me_g(\Q)$ for all $g$ 
up to 13.

\begin{cor}
  For each $g \le 13$, the curves of genus $g$ over $\Q$ for which 
  the group generated by the Weierstrass points in its Jacobian is 
  isomorphic to $\Z^{g(g^2-1)-1}$, determine a Zariski dense set of 
  moduli in $\Me_g$.
\end{cor}

\begin{proof}
For each $g \le 13$ the moduli space $\Me_g$ is 
unirational~\cite{ArbarelloSernesi79,ChangRan84,Sernesi81}, 
i.e.~$\Me_g$ is covered by a dominant map $\pi:\PP^N \rightarrow \Me_g$. 
Let $\Ae \rightarrow \Me_g^\circ$ be as above, and let $\Ae_U = 
\Ae \times_\pi U$ be the base extension to $U = \pi^{-1}(\Me_g^\circ)$. 
The Weierstrass subgroup embeds in $\Ae_U(\Q(\PP^N)) = \Ae(\Q(\PP^N))$,
which specialises injectively to $(\Ae_U)_P(\Q) = \Ae_{\pi(P)}(\Q)$ for 
all $P$ outside of a thin set $Z$ in $U(\Q)$.  The image $\pi(U(\Q)-Z)$ 
is then a Zariski dense set in $\Me_g$ consisting of moduli of curves 
for which the Weierstrass subgroup attains the maximal rank.
\end{proof}

Lang conjectures that for a variety $X$ of general type over a number 
field $F$, for each finite extension $K/F$, the Zariski closure of $X(K)$ 
is contained in a proper subvariety of $X$ (see~\cite[pp. 15--20]{Lang91}
or~\cite[\S F.5.2]{HS00}).
Since the variety $\Me_g$ is known to be of general type for $g \ge 24$ 
(see \cite{EisenbudHarris87b,Harris84,HarrisMumford82}), Lang's conjecture 
would imply that for any number field $K$ the set $\Me_g(K)$ is contained 
in a proper closed subvariety of $\Me_g$.  Since the complete set of 
rational points over $K$ is not dense, the density result can not be 
expected to hold for large $g$. 
In contrast, the plane curves of degree $d$, determine unirational 
subvarieties of $\Me_g$, for $g = (d-1)(d-2)/2$.  An analogous result to 
Theorem~\ref{thm_monodromy} for the Galois group of Weierstrass points 
of plane curves of degree $d$ (in the spirit of~\cite{Harris79}) would 
be desirable in order to establish the maximality of the rank of the 
Weierstrass subgroup for plane curves.
\vspace{4mm}

\noindent{\it Acknowledgements.}  The authors thank Marc Hindry and Ren{\'e}
Schoof for interest and comments on an earlier draft of this work.

\vspace{8mm}
\noindent
\begin{minipage}[t]{62mm}
{\it 
Martine Girard, David Kohel \\
School of Mathematics and Statistics \\
The University of Sydney NSW 2006 \\
Australia}\\
{\tt girard@maths.usyd.edu.au},\\
\quad{\tt kohel@maths.usyd.edu.au}
\end{minipage}
\begin{minipage}[t]{62mm}
{\it 
Christophe Ritzenthaler\\
Institut de Math\'ematiques de Luminy \\
163 Avenue de Luminy, Case 907\\
13288 Marseille Cedex 9\\
France}\\
{\tt ritzenth@math.jussieu.fr}
\end{minipage}

\end{document}